\documentclass[12pt,a4paper]{article}
\usepackage{amsmath,amsfonts,makeidx}
\makeindex
\numberwithin{equation}{section}
\newtheorem{thm}{Theorem}[section]
\newtheorem{cor}[thm]{Corollary}
\newtheorem{lem}[thm]{Lemma}
\newtheorem{prop}[thm]{Proposition}
\newtheorem{rem}[thm]{Remark} 
 
\newtheorem{defn}[thm]{Definition} 
 
\newtheorem{conj}[thm]{Conjecture} 
\newtheorem{hyp}[thm]{Hypothesis}

\makeatletter

\def\section{\@startsection {section}{1}{\z@}{3.5ex plus 1ex minus
    .2ex}{2.3ex plus .2ex}{\large\bf}}
    \def\subsection{\@startsection{subsection}{2}{\z@}{3.25ex plus 1ex minus 
 .2ex}{1.5ex plus .2ex}{\bf}}

\makeatother

\newcommand{\qed}{\hspace{2em}\nolinebreak\rule{1.5mm}{3.5mm}}

\def\S{\mathfrak{S}}
\def\Q{\mathbb{Q}}
\def\N{\mathbb{N}}

\def\C{\mathbb{C}}

\def\Z{\mathbb{Z}}

\def\ker{{\rm Ker}}
\def\res{{\rm Res}}
\def\reg{{\text{\rm reg}}}
\def\deg{{\rm deg}}
\def\diag{{\rm diag}}
\def\id{{\rm Id}}
\def\tr{{\rm tr}}
\def\class{{\rm Class}}
\def\rep{{\rm Rep}}
\def\harm{{\rm Harm}}
\def\Hom{{\rm Hom}}
\def\funct{{\rm Funct}}

\begin{document}

\pagestyle{plain} 
\title{Complex Reflection Groups and Fake Degrees\thanks{Lectures at
RIMS, Kyoto University (Japan) in 1997. Noted by Kenji Taniguchi.}}

\author{Eric M. Opdam\thanks{The author thanks Gunter Malle 
and Kenji Taniguchi for their keen observations, questions and comments, 
and Kenji Taniguchi for his commitment and skill in writing these 
notes.}}

\date{version 98.4.29}

\maketitle
\tableofcontents

\section{Introduction}

{\small In this section, we introduce finite complex reflection 
groups, define the fake degrees, and explain the objectives of 
this paper.}
\vskip3mm

\begin{defn} 
{\rm Let  \index{v@$V$} $V$ be an $n$-dimensional Hilbert space 
and let \index{uv@$U(V)$} $U(V)$ be the group of unitary transformations on $V$. 
An element $r\in U(V)$ is called a \index{complex reflection} 
{\it complex reflection} 
if the set $\ker(r-\id_V)$ of fixed points 
of $r$ in $V$ is a complex hyperplane $H_r$. 
A finite subgroup $W$ of $U(V)$ is called 
\index{finite complex (or unitary) reflection group} 
{\it a finite complex 
(or also unitary) reflection group} if $W$ is generated by 
complex reflections.\rm} 
\end{defn}

The following theorem is the fundamental fact characterizing 
such subgroups in $U(V)$. 

\begin{thm}[Shephard-Todd \cite{ST}]\label{thm:sh-to}
A subgroup $W$ of $U(V)$ is a finite complex reflection group 
if and only if the subring \index{pw@$P^W$} 
$P^W$ of the $W$-invariant elements 
in the ring \index{p@$P$} 
$P$ of polynomial functions on $V$ is generated by $n$ 
algebraically independent homogeneous elements \index{pi@$p_i$} 
$p_1,\dots,p_n$. 
\end{thm}

In this case, the homogeneous degrees \index{di@$d_i$} 
$d_i=\deg(p_i)$ depend only on $W$,
and are called the \index{primitive degree} 
{\it primitive degrees} of $W$. 
The product of $d_i$ is equal to $|W|$ and the sum 
$\sum_{i=1}^n(d_i-1)$ is equal to the number of complex reflections 
in $W$. 

We call $W$ {\it irreducible} if $W$ acts irreducibly on $V$.
Shephard and Todd have given the complete classification of the 
irreducible complex reflection groups:

\begin{thm}[Shephard-Todd \cite{ST}] If $W$ is an irreducible complex 
reflection group, then $W$ is isomorphic to one of the following list:
{\rm \begin{enumerate} 
\item
The symmetric Group $\S_n$, acting on $V=\{v\in \C^n\mid \sum v_i=0\}$.
 
\item
Let $m, p, n$ be positive integers such that $p$ divides $m$, 
$m\geq 2$, and $p=1$ if $n=1$. 
Let $\{e_i; i=1,\dots,n\}$ be an orthonormal basis of $\C^n$ 
and for each $i=1,\dots,n$, let $\zeta_i$ be a $m$-th 
root of unity,  such that $\left(\prod_{i=1}^n\zeta_i\right)^{\frac{m}{p}}=1$. 
We denote by \index{gmpn@$G(m,p,n)$} $G(m,p,n)$ the group generated by 
\[e_i\mapsto \zeta_ie_{\sigma(i)}\quad(\sigma\in\S_n).\]
This is a finite complex reflection group in $\C^n$. These groups are
the 
\index{imprimitive complex reflection group} 
{\it imprimitive} complex reflection groups.
 
\item
One of 34 exceptional cases (these cases contain of course the 
exceptional Coxeter groups). 
\end{enumerate}}
\end{thm} 

\begin{rem} 
{\rm Many rank two cases are obtained as follows: 
Let $\Gamma\subset SL_2(\C)$ be a finite subgroup 
corresponding to a platonic polyhedron. 
For suitable choice of $a\in\N$, 
$\Gamma\cdot \mu_a$ is a finite complex reflection group, 
where $\mu_a=\left\{\pmatrix\alpha&0\\0&\alpha^{-1}\endpmatrix; 
\alpha^a=1\right\}$.} 
\end{rem} 
\vskip3mm

The action of $W$ on $P$ respects the grading by degree. 
Let \index{pp@$P_+$} $P_+$ 
be the graded ideal of polynomials vanishing at $0\in V$ and 
let \index{ppw@$P_+^W$} $P_+^W$ 
be the ideal generated by the invariants in $P_+$. 
The \index{coinvariant algebra} coinvariant algebra 
$P/PP_+^W$ is a representation space of $W$ 
and we denote it by \index{rw@$RW$} $RW$. 
This representation is isomorphic to 
the left regular representation on $\C[W]$. 
The space $RW$ is graded by homogeneous degree and this grading is 
compatible with the action of $W$. 
Let \index{rwk@$RW_k$} $RW_k$ be the homogeneous subspace of $RW$ of degree $k$. 

It is well known that the \index{graded character} 
graded character of this representation is 
\index{trw@$\tr_{RW}(w)$} 
\begin{equation} 
\tr_{RW}(w):=\sum_{k\geq0}(\tr(w)|_{RW_k})T^k
=\frac{\prod_{i=1}^n(1-T^{d_i})}{\det_V(1-Tw^{-1})}
\in\C[T].\label{eq:graded char}
\end{equation}
By this character, we have a map from the space $\class(W)$ 
of class functions on $W$ 
to $\C[T]$: 
\[\class(W)\ni\alpha\mapsto 
F_\alpha:=\frac1{|W|}\sum_{w\in W}\alpha(w)\tr_{RW}(w^{-1})\in\C[T].\]

\begin{defn}\label{fake}
{\rm For a representation $\tau$ of $W$ we write 
\index{rt@$R_\tau$} \index{ft@$F_\tau$} 
$R_\tau=F_{\overline \tau}=F_{\overline{\chi_\tau}}$, where 
\index{chit@$\chi_\tau$} $\chi_\tau$ 
denotes the character of $\tau$. $R_\tau$ is called the 
\index{fake degree} 
{\it fake degree} of $\tau$. Notice that $R_\tau\in\Z[T]$ with nonnegative 
coefficients.}
\end{defn}

Since $RW\simeq\C[W]$, $\dim\Hom_W(\tau,RW)=\deg\tau=:l$ for 
each irreducible representation $\tau\in\hat{W}$. 
Let \index{pti@$p_i^\tau$} 
$p_1^\tau\leq\dots\leq p_l^\tau$ be the homogeneous embedding 
degrees of $\tau$ in the coinvariant algebra $RW$, i.e. 
the degrees of irreducible $\tau$-components in $RW$.

\begin{cor}\label{cor:1.6}
For an irreducible representation $\tau\in\hat{W}$, 
\begin{align}
F_\tau(T)=&\frac1{|W|}\sum_{w\in W}\frac{\chi_\tau(w)}
{\det_V(1-Tw)}\prod_{i=1}^n(1-T^{d_i})\label{eq:1.2}\\
=&\sum_{j=1}^lT^{p_j^\tau}.\label{eq:1.3}
\end{align}
\end{cor}
\begin{pf} 
By \eqref{eq:graded char} and by the orthogonality of characters. 
\qed
\end{pf}
 
The name {\it fake degree} was given by Lusztig. When $W$ is a Coxeter group the  
notion of fake degree plays a role 
in the representation theory of \index{finite Chevalley group} 
finite Chevalley groups. By Lusztig's work,
they can be regarded as approximations of the degrees of the 
principal series \index{unipotent representation} 
unipotent representations of a finite Chevalley group. In fact 
Lusztig has shown that every \index{unipotent degree} 
unipotent degree of a finite Chevalley 
group \index{gfq@$G({\mathbb F}_q)$} $G({\mathbb F}_q)$ 
with Weyl group $W$, 
can be expressed as a rational linear combination 
of the fake degrees of $W$ evaluated at $T=q$ (see \cite{C}, \cite{L}). 

The ideas and conjectures of Brou\'e, Malle and Michel 
state, among many other things, that this role of $W$ and its fake degrees 
for the study of unipotent representations can be extended naturally to the 
more general case where $W$ is a complex reflection group that arises as the 
quotient $N(L)/L$ with $L\subset G$ a  
\index{d-cuspidal Levi subgroup} 
``d-cuspidal'' Levi subgroup of $G$. 
We refer the reader to \cite{BM}, \cite{BMM} for an account of 
this subject. Also see \cite{Ma}, where the unipotent degrees and fake degrees 
were studied in the case of the general imprimitive group.

A central notion in these considerations of unipotent and fake degrees for 
a complex reflection group $W$  is its 
\index{cyclotomic Hecke algebra} {\it (cyclotomic) Hecke algebra}. 
It is a deformation of the group algebra of $W$, 
similar to the ordinary Hecke algebra 
of a Coxeter group. It was introduced by Brou\' e and Malle in  
\cite{BM}, and investigated in many subsequent papers (for example 
\cite{Ma}, \cite{M}, \cite{BMR}). From the results of these papers it is clear that the cyclotomic 
Hecke algebra shares many of the properties  
which give the usual Hecke algebra
its prominent role in representation theory. But also, there are still 
severe problems to give the cyclotomic Hecke algebra a transparent theoretical 
basis similar to the theory of ordinary Hecke algebras. Many results rely 
on classifications and computer aided computations. In fact there are many exceptional 
complex reflection groups for which it is still not known whether their 
Hecke algebras are free over the 
coefficient ring or not, or what the rank of the Hecke algebra is. 
The situation is better for the 
imprimitive cases. Here the theory is rather well understood, see \cite{AK} and 
\cite{Ma}. 

In these lectures we approach the cyclotomic Hecke algebra from topology. 
It is the same approach as was used in the paper \cite{BMR}.  
This way of thinking 
about the Hecke algebra is quite natural in the case of a Coxeter group $W$, given 
Brieskorn's description of the fundamental group of the regular orbit space 
as the braid group of $W$ \cite{Br}. 
For Coxeter groups one can find results in this 
direction in  \cite{HO}, \cite{K1} (and many other papers).   
We use the monodromy representation of certain systems of differential equations 
to construct the ``topological cyclotomic Hecke algebra'' (terminology from 
\cite{BMR}), which is known to be isomorphic to the cyclotomic Hecke algbra 
in many cases (and conjectured to be isomorphic in all cases) (see \cite{BMR}).

There are two main results in these notes. 
The first is described in section 4. This is a transformation 
property of fake degrees of representations of $W$ with respect to  
certain operations on representations of $W$.
In view of the current difficulties with cyclotomic Hecke algebras we  
formulate and prove this result without the use of cyclotomic Hecke algebras. 
The price we have to pay is that the ``operations'' on representations of $W$ are 
somewhat mysterious without the cyclotomic Hecke algebra, and not much can be said about 
their basic properties. In the second half of these notes we assume certain 
facts about the cyclotomic Hecke algebra, and with these assumptions we 
interpret the ``operations'' mentioned above as the natural action on ${\rm Irr}(W)$ 
of the \index{geometric Galois group} 
geometric Galois group of the character field of the cyclotomic Hecke algebra. 
See section 7 for more details about this part of the story. The main result of this 
second part is Theorem~\ref{thm:pol}, which implies that the geometric 
Galois group of the character field over the coefficient field is 
abelian (but the statement of Theorem~\ref{thm:pol} is more precise).



\section{Minimal polynomial realization} 

{\small In this section, we define the minimal $\tau$-matrix 
and investigate the basic properties of it.} 
\vskip3mm

Let \index{c@$\mathcal{C}$} 
$\mathcal{C}$ be the set of $W$-orbits of reflection hyperplanes 
and let \index{a@$\mathcal{A}$} 
$\mathcal{A}$ be the full hyperplane arrangement 
$\cup_{C\in\mathcal{C}}C$. 
For a hyperplane $H$ in an orbit $C\in{\mathcal{C}}$, 
the stabilizer \index{wh@$W_H$} 
\[W_H:=\{w\in W; wx=x\enskip\text{for any}\enskip x\in H\}\]
is isomorphic to a cyclic group $\Z/e_C\Z$. 
Here, \index{ec@$e_C$} 
$e_C$ is the order of $W_H$, which is determined by the orbit $C$. 
For a representation $\tau$ of $W$, 
we define a nonnegative integer \index{ncjt@$n_{C,j}^\tau$} $n_{C,j}^\tau$ by 
\begin{equation}\label{eq:rest}
\res_{W_H}^W\tau\simeq\bigoplus_{j=0}^{e_C-1}n_{C,j}^\tau\det{}^{-j}.
\end{equation}
Note that $\hat{W}_H=\{1=\det^{-e_C}, \det^{-e_C+1},\dots,\det^{-1}\}$. 

In this note, we denote $f(H)=f(C)$ for a function 
$f$ on $\mathcal{C}$ if $H\in C\in{\mathcal{C}}$. 
For example, $e_H=e_C$ and $n_{H,j}^\tau=n_{C,j}^\tau$. 
Note that 
\begin{align}
\sum_{i=1}^n(d_i-1)&=|\{\text{\rm complex reflections}\}|
=\sum_{H\in{\mathcal{A}}}(e_H-1),\label{eq:refno}\\
\deg\tau&=\sum_{j=0}^{e_H-1}n_{H,j}^\tau
\quad(\text{\rm for each $H\in{\mathcal{A}}$}).\label{eq:deg-n}
\end{align}

\begin{lem}\label{lem:p-n}
For an irreducible representation $\tau$, 
\begin{equation}
\sum_{j=1}^lp_j^\tau
=\sum_{C\in{\mathcal{C}}}|C|\sum_{j=1}^{e_C-1}jn_{C,j}^\tau.
\label{eq:p-n}
\end{equation}
\end{lem}

\begin{pf}
To show this, we use Corollary~\ref{cor:1.6}. 
We differentiate this equality by $T$ and take the limit $T\rightarrow1$. 
After we take the limit,
\eqref{eq:1.3} is equal to $\sum_{j=1}^lp_j^\tau$, 
and the only terms which survive 
in the sum $\sum_{w\in W}$ of right hand side of \eqref{eq:1.2} are 
$w=e$ and $w=$complex reflections. 
The limit of $w=e$ term is 
\begin{align*}
\lim_{T\rightarrow1}\frac{\deg\tau}{|W|}&
\sum_{i=1}^n\left(\sum_{j_i=0}^{d_i-1}j_iT^{j_i-1}\right)
\prod_{k\not=i}\left(\sum_{j_k=0}^{d_k-1}T^{j_k}\right)\\
=&\frac{\deg\tau}{|W|}\sum_{i=1}^n\frac12d_i(d_i-1)
\prod_{k\not=i}d_k
=\sum_{H\in{\mathcal{A}}}\sum_{j=0}^{e_H-1}\frac{e_H-1}2n_{H,j}^\tau.
\end{align*}
Here, we used the facts \eqref{eq:refno}, \eqref{eq:deg-n} 
and $\prod_{i=1}^nd_i=|W|$. 

For $H\in{\mathcal{A}}$, we denote by \index{zeh@$\zeta_H$} 
$\zeta_H$ a primitive 
$e_H$-th root of unity. 
The limit of the sum of $w=$complex reflections is 
\begin{align*}
\lim_{T\rightarrow1}&\frac{d}{dT}\frac1{|W|}
\sum_{H\in{\mathcal{A}}}
\sum_{k=1}^{e_H-1}\sum_{j=0}^{e_H-1}
\frac{n_{H,j}^\tau\zeta_H^{-jk}}{1-T\zeta_H^k}
(1-T)\prod_{i=1}^n\left(\sum_{j_i=0}^{d_i-1}T^{j_i}\right)\\
&=-\sum_{H\in{\mathcal{A}}}
\sum_{k=1}^{e_H-1}\sum_{j=0}^{e_H-1}
\frac{n_{H,j}^\tau\zeta_H^{-jk}}{1-\zeta_H^k}\\
&=\dots\\
&=-\sum_{H\in{\mathcal{A}}}\sum_{j=0}^{e_H-1}n_{H,j}^\tau
\frac{e_H-1-2j}2.
\end{align*}
\qed
\end{pf}

\begin{defn}\label{defn:tau-m}
{\rm For a finite dimensional representation $(\tau, E)$ of $W$, 
we choose an explicit matrix realization $\tau: W\rightarrow GL(l,\C)$ 
($l=\deg\tau$) and fix it. 

A matrix \index{m@$M$} $M=(m_{ij})\in Mat(l\times l,P)$ with $\det(M)\not=0$ 
is called a \index{minimal $\tau$-matrix} {\it minimal $\tau$-matrix} 
if it satisfies the following conditions: 
\begin{enumerate}
\item
For each $w\in W$, $M^w:=(m_{ij}\circ w^{-1})=\tau(w)M$. 
\item
Let \index{e@${\mathbb E}$} 
${\mathbb E}:=\sum_{i=1}^nx_i\partial/\partial x_i$ be the 
Euler vector field. The matrix $M$ satisfies 
\[{\mathbb E}M:=({\mathbb E}m_{ij})=M\cdot C_{\mathbb E}^M\]
for some \index{cme@$C_{\mathbb E}^M$} 
$C_{\mathbb E}^M\in Mat(l\times l,\C)$. 
\item
\[\tr(C_{\mathbb E}^M)=\sum_{C\in{\mathcal{C}}}|C|\sum_{j=1}^{e_C-1}jn_{C,j}^\tau.\]
\end{enumerate}}
\end{defn}

\begin{rem}\label{rem:mintau}
{\rm 
\begin{enumerate}
\item \ref{defn:tau-m} Condition (iii) is equivalent to the following (iii)$'$:
\[{\rm deg}(\det(M))=\sum_{C\in{\mathcal{C}}}|C|\sum_{j=1}^{e_C-1}jn_{C,j}^\tau.\]
This follows immediately by considering the action of the group
$\C^\times$ (for 
which ${\mathbb E}$ is the infinitesimal generator) on $M$.
\item If $M$ is a minimal $\tau$-matrix and $g\in GL(l,\C)$, 
then $Mg$ is also a minimal $\tau$-matrix.
\end{enumerate} 
}
\end{rem}

\begin{prop}\label{prop:nature of mintau}  
\begin{enumerate}
\item
$C_{\mathbb E}^M$ is semisimple and its spectrum is contained in the set of 
non-negative integers. 
\item
For every finite dimensional representation $\tau$ of $W$, 
there exists a minimal $\tau$-matrix $M$. 
\item
Let \index{alh@$\alpha_H$} $\alpha_H$ be a linear function 
satisfying $\ker\alpha_H=H$. 
We define \index{pic@$\pi_C$} 
$\pi_C=\prod_{H\in C}\alpha_H$ for $C\in{\mathcal{C}}$. 
Then 
\begin{equation}
\det(M)=const.\prod_{C\in{\mathcal{C}}}\pi_C^{\sum_{j=0}^{e_C-1}jn_{C,j}^\tau}.
\label{eq:detM}\end{equation}
\item
For any $N\in Mat(l\times l, P)$ satisfying Definition~\ref{defn:tau-m} (i), 
there exists $R\in\
 Mat(l\times l,P^W)$ such that $N=MR$. 
\item
Spectrum of $C_{\mathbb E}^M$ does not depend on $M$, only on $\tau$.
\end{enumerate}
\end{prop} 

\begin{pf} 
(i) Since $e^{2\pi\sqrt{-1}{\mathbb E}}$ acts on $P$ by identity, we have 
$e^{2\pi\sqrt{-1}C_{\mathbb E}^M}=\id$. 
It follows that $C_{\mathbb E}^M$ is semisimple and 
its eigenvalues are integers. 
The spectrum of $C_{\mathbb E}^M$ is contained 
in the set of non-negative integers 
since $M$ has polynomial entries. 

(ii) For the proof of (ii), we may assume $\tau$ to be irreducible. 
Let \index{harmt@$\harm(\tau)$} 
$\harm(\tau)$ be the vector space spanned by harmonic polynomials 
associated with $\tau$ and let \index{taus@$\tau^*$} 
$(\tau^*, E^*)$ be the contragredient 
representation of $(\tau, E)$. Take bases \index{ei@$\varepsilon_i$} 
$\varepsilon_i$ of $E$ and \index{sj@$\sigma_j$} 
$\sigma_j$ of $(\harm(\tau)\otimes E^*)^W$ ($i,j=1,\dots,l$) and 
we define \index{hijt@$h_{ij}^\tau$} 
$M=(h_{ij}^\tau):=(\sigma_j(\varepsilon_i))$. 
By this construction, (i) and (ii) in Definition~\ref{defn:tau-m} are clear. 
We have to check Definition~\ref{defn:tau-m} (iii), 
but this is clear by 
$\tr(C_{\mathbb E}^M)=\sum_{j=1}^l\deg h_{jj}^\tau
=\sum_{j=1}^lp_j^\tau$ (using Lemma \ref{lem:p-n}). 

(iii) By \eqref{eq:rest} and Definition \ref{defn:tau-m} (i), 
there exists $D\in GL(l,\C)$ such that \index{mh@$M_H$} 
$D\tau(s)D^{-1}$ ($s\in W_H$) is diagonal and 
\[DM=M_HM',\]
where $M'$ is a $W_H$-invariant matrix with polynomial entries and 
\begin{equation} 
M_H=\diag(I_{n_{H,0}^\tau},
\alpha I_{n_{H,1}^\tau},
\dots,\alpha_H^{e_H-1}I_{n_{H,e_H-1}^\tau}).
\end{equation}
($I_p$ is the $p$-th unit matrix.) 
Note that $M_H$ is a minimal $\tau|_{W_H}$-matrix and
$\det(M')\not=0$. 

By the above discussion, $\det(M)$ is divisible by 
$\alpha_H^{\sum_{j=0}^{e_H-1}jn_{H,j}^\tau}$ 
for each $H\in\mathcal{A}$ and also by the right hand side of \eqref{eq:detM}. 
Both hand sides of \eqref{eq:detM} have the same degree since 
$\deg(\det(M))=\tr(C_{\mathbb E}^M)=\sum_{j=1}^lp_j^\tau
=\sum_{C\in{\mathcal{C}}}|C|\sum_{j=1}^{e_C-1}jn_{C,j}^\tau
=\deg\left(\prod_{C\in{\mathcal{C}}}
\pi_C^{\sum_{j=0}^{e_C-1}jn_{C,j}^\tau}\right)$. 
This proves \eqref{eq:detM}. 

(iv) By the same discussion as above, 
we can show that, for any $N\in Mat(l\times l,P)$ satisfying (1), 
there exists a $W_H$-invariant matrix $N'\in Mat(l\times l, P)$ 
such that $DN=M_HN'$. 
Hence, for each reflection hyperplane $H$, 
$M^{-1}N=(M')^{-1}N'$ is $W_H$-invariant and regular at $H$. 
It follows that $M^{-1}N$ is $W$-invariant and the entries of $M^{-1}N$ 
are polynomials because of \eqref{eq:detM}. 

(v) By Remark \ref{rem:mintau} and (i) of this proposition, 
we may assume that $C_{\mathbb{E}}^{M_\nu}$ ($\nu=1,2$) 
are diagonal matrices. 

Let $n_1^\nu, \dots, n_l^\nu$ be the spectrum of $C_{\mathbb E}^{M_\nu}$ 
($\nu=1,2$). 
Since 
\[a^{C_{\mathbb E}^{M_\nu}}=\diag(a^{n_1^\nu},\dots,a^{n_l^\nu}),\]
the matrix $R=M_1^{-1}M_2$ satisfies $R(ax)_{ij}=a^{n_j^2-n_i^1}R(x)_{ij}$,  
and it follows that $R(x)_{ij}=0$ if $n_j^2<n_i^1$. 
But we know $\sum_{i=1}^ln_i^1=\sum_{i=1}^ln_i^2$ by 
Definition~\ref{defn:tau-m} (iii). 
If the sequences $(n_1^\nu,\dots,n_l^\nu)$ ($\nu=1,2$) do not coincide, 
then, for every $\sigma\in\S_l$, there exists $i$ such that 
$n_{\sigma(i)}^2<n_i^1$. 
This implies $\det R=0$, which is a contradiction. \qed
\end{pf}

\begin{rem} 
{\rm This generalizes a construction of Stanley \cite{S}, who proved 
that for every one dimensional representation $\tau\in\hat W$ the  
\index{pseudo-invariant} pseudo-invariants of type $\tau$ 
in $P$ form a rank one free module over $P^W$, with generator 
\[\prod_{C\in{\mathcal{C}}}\pi_C^{\sum_{j=0}^{e_C-1}jn_{C,j}^\tau}\]
which is the minimal $\tau$-matrix in this situation. So Stanley's result 
is Proposition \ref{prop:nature of mintau} (iv) in this special case.}
\end{rem}

\begin{cor}\label{cor:2.6}
Let $M$ be a minimal $\tau$-matrix and let \index{ni@$n_i$} 
$n_1\leq\dots\leq n_l$ 
be the spectrum of $C_{\mathbb E}^M$. 
Then 
\[F_\tau(T)=\sum_{i=1}^lT^{n_i}.\]
\end{cor}



\section{Knizhnik-Zamolodchikov equations}

{\small 
We are going to construct deformed minimal $\tau$-matrices using certain 
differential equations. 
This is a generalization of the construction of a minimal $\tau$-matrix 
using harmonic polynomials.}  
\vskip3mm

Let us choose labels 
\index{k@$k=(k_{C,j})_{C\in{\mathcal{C}}, j=0,\dots,e_C-1}$} 
$k=(k_{C,j})_{C\in{\mathcal{C}}, j=0,\dots,e_C-1}$ 
with $k_{C,j}\in\C$. 
Put \index{qcj@$q_{C,j}=\exp(-2\pi\sqrt{-1}k_{C,j})$} 
$q_{C,j}=\exp(-2\pi\sqrt{-1}k_{C,j})$. We shall sometimes use the 
notation \index{q@$q=(q_{C,j})_{C\in{\mathcal{C}}, j=0,\dots,e_C-1}$} 
$q$ for the vector 
$(q_{C,j})_{C\in{\mathcal{C}}, j=0,\dots,e_C-1}$. 

Let \index{ejh@$\varepsilon_j(H)$} 
$\varepsilon_j(H)$ be the idempotent element 
$\frac1{e_H}\sum_{w\in W_H}\det^j(w)w$ in $\C[W_H]$. 

\begin{defn} 
{\rm We define \index{omega@$\omega$} 
\[\omega=\sum_{H\in{\mathcal{A}}}a_H\omega_H,\]
where \index{ah@$a_H$} 
$a_H=\sum_{j=0}^{e_H-1}e_Hk_{H,j}\varepsilon_j(H)$ and 
\index{omegah@$\omega_H$} 
$\omega_H=d(\log\alpha_H)=\frac{d\alpha_H}{\alpha_H}$. }
\end{defn}

Let \index{vreg@$V^\reg$} 
\[V^\reg=\{v\in V; v\not\in H\enskip \text{for any}\enskip H\in\mathcal{A}\}\] 
and let us denote by \index{ow@$\mathcal{O}[W]$} 
$\mathcal{O}[W]$ (resp. \index{omegaw@$\Omega^1[W]$} $\Omega^1[W]$) 
the sheaf of germs of $\C[W]$-valued holomorphic functions (resp. 
holomorphic $1$-forms) on $V^\reg$. 
Then $\omega$ is an element of $\Omega^1(V^{\reg})[W]$ satisfying 
$w\cdot(\omega\circ w^{-1})\cdot w^{-1}=\omega$.


\begin{thm}[Kohno, Brou\'e-Malle-Rouquier, Opdam]\label{thm:kbmro}
\begin{enumerate}
\item
$\omega$ is integrable one form, i.e. $\omega\wedge\omega=0$. 
This means that the connection \index{nablak@$\nabla(k)$} 
\[\nabla(k): {\mathcal{O}}[W]\ni\Phi\mapsto d\Phi+\omega\Phi
\in\Omega^1[W]\]
is completely integrable, i.e. the $\nabla(k)$-flat local sections 
$\Phi\in{\mathcal{O}}[W]$ form a vector space of $\dim|W|$. 
\item
The connection $\nabla(k)$ commutes with right $W$-multiplication. 
\item
We define an action $\Phi\mapsto\Phi^w$ of $W$ on 
${\mathcal{O}}[W]\simeq{\mathcal{O}}\otimes\C[W]$ and
$\Omega^1[W]\simeq\Omega^1\otimes\C[W]$ by 
$w\otimes(\text{\rm left multiplication of $w$})$. 
Then $\nabla(k)$ commutes with this action. 
\end{enumerate}
\end{thm} 

\begin{pf} 
Basically, the proof of Kohno \cite{K} is still valid.\qed 
\end{pf} 

Let $(\tau, E)$ be a representation of $W$,  
and let \index{et@$e_\tau$} $e_\tau$ be an idempotent of the ring 
$\C[W]^N$ such that $\C[W]^N\cdot e_\tau$ is isomorphic to 
$(\tau, E)$ (we choose $N$ large enough). 
By Theorem \ref{thm:kbmro} (ii), 
$\nabla(k)^N$ descends to a connection on the bundle 
$V^\reg\times E\rightarrow V^\reg$. 
Moreover, by Theorem \ref{thm:kbmro} (iii), $\nabla(k)^N$ descends to 
\index{e@${\mathcal{E}}$} 
${\mathcal{E}}:={\mathcal{O}}(X^\reg)\otimes_\C{\mathcal{L}}(E)$, 
where \index{x@$X$} $X=W\backslash V$, \index{xreg@$X^\reg$} 
$X^\reg=W\backslash V^\reg$ and \index{le@${\mathcal{L}}(E)$} 
${\mathcal{L}}(E)$ is the local system $V^\reg\times_WE\rightarrow
X^\reg$. The resulting connection on ${\mathcal{E}}$, 
called \index{Knizhnik-Zamolodchikov (KZ) connection} 
{\it Knizhnik-Zamolodchikov (KZ) connection}, 
is denoted by 
\index{nablatk@$\nabla_\tau(k)$} $\nabla_\tau(k)$.

\begin{defn}\label{kztau}
{\rm We define \index{enabla@${\mathcal{E}}^{\nabla_\tau(k)}$} 
${\mathcal{E}}^{\nabla_\tau(k)}\rightarrow X^\reg$ to be the local 
system of $\nabla_\tau(k)$-flat sections in $\mathcal{E}$.} 
\end{defn}

Let \index{xzero@$x_0$} $x_0$ be a base point in $X^\reg$, 
and \index{vzero@$v_0$} $v_0$ a lift of $x_0$ 
in $V^\reg$. 
For every $H\in{\mathcal{A}}$, 
we choose a path \index{lh@$l_H$} $l_H: v_0\rightarrow s_Hv_0$ in $V^\reg$, 
where \index{sh@$s_H$} $s_H$ is a generator for $W_H$ such that 
$\det(s_H)=\zeta_H=e^{2\pi\sqrt{-1}/e_H}$. 
We denote $\pi_1(V^\reg,v_0)$ and $\pi_1(X^\reg,x_0)$ by \index{p@$P$}
$P$ and \index{b@$B$} $B$, 
and we call them the \index{pure braid group} 
{\it pure braid group} and the \index{braid group} {\it braid group}, 
respectively.

The following theorem is due to Brou\'e, Malle and Rouquier: 

\begin{thm}[see \cite{BMR}]\label{thm:bmr}
\begin{enumerate} 
\item
$B$ is generated by $\{l_H\}_{H\in{\mathcal{A}}}$. 
\item
$P$ is generated by $\{l_H^{e_H}\}_{H\in{\mathcal{A}}}$.
\item
We have the following short exact sequence:
\[1\rightarrow P\rightarrow B\rightarrow W\rightarrow1,\]
where the map $B\rightarrow W$ is given by $l_H\mapsto s_H$. 
\end{enumerate}
\end{thm}

\begin{defn}\label{defn:monodromy} 
{\rm Let \index{tauk@$\tau(k)$} $\tau(k)$ be the 
monodromy action of $B$ on 
${\mathcal{E}}^{\nabla_\tau(k)}_{x_0}$. }
\end{defn}

\begin{thm}[Brou\'e-Malle-Rouquier, Opdam]\label{thm:bmro} 
\[\prod_{j=0}^{e_H-1}(\tau(k)(l_H)-q_{H,j}\zeta_H^j)=0.\]
(Note that 
$q_{H, j}\zeta_H^j=\exp\left(2\pi\sqrt{-1}(j-e_Hk_{H,j})/e_H\right)$.)
Moreover, with respect to a fixed basis of ${\mathcal{E}}_{x_0}=E$
 the matrix of 
$\tau$ will have coefficients in \index{s@$S$} $S$, the ring of 
entire functions in the labels 
$(k_{C,j})_{C\in{\mathcal{C}}, j=0,\dots,e_C-1}$.
\end{thm}

\begin{pf}
The proof we present here differs a little bit from the one in \cite{BMR}, 
and will give important additional information about the local 
behaviour of flat sections near the reflection hyperplanes.

We first address the last assertion of Theorem~\ref{thm:bmro}. 
This is a basic fact, which we prove anyway, for want of a good 
reference. 
The connection matrix of $\nabla_\tau$ depends polynomially on the labels 
$(k_{C,j})_{C\in{\mathcal{C}}, j=0,\dots,e_C-1}$. By definition of 
analytic continuation along a path, it therefore suffices to prove the 
following local fact (the monodromy matrix of a loop is a composition 
of finitely many local steps like these):   
Suppose we have a holomorphic, first order, linear system of 
differential equations on the bundle 
$\C^l\times D$ over the 
unit disc $D$, and the coefficient matrix depends polynomially on 
a parameter $\kappa\in\C$. 
Given $v\in \C^l$, the unique solution $\nu(z,\kappa)$ ($z\in D$) such 
that $\nu(0,\kappa)=v$ is an entire function of $\kappa$. 
For this it suffices to check that the power series expansion of $\nu$ 
on $D$ converges locally uniformly in $\kappa$, and this is an elementary 
exercise left to the reader. 

Choose $H_0\in{\mathcal{A}}$ and fix it. For notational convenience, 
we abbreviate $\alpha_{H_0}$ as $\alpha_0$, $s_{H_0}$ as $s_0$ and so on. 
Let $x_0\in H_0$ be a regular point and let $(x,\alpha_0)$ be a coordinate 
in a tubular neighborhood $U\times I$ of $x_0$, 
where $U\subset H_0$ and $I\subset\C\alpha_0$. 
For every $\varepsilon_0\in E$ with 
$s_0\varepsilon_0=\zeta_0^{-j}\varepsilon_0$, 
we shall construct a flat section 
\index{exa@$\varepsilon(x,\alpha_0)$} 
$\varepsilon(x,\alpha_0)$ in $U\times I$. 

Contracting the equation $d\varepsilon+\omega\varepsilon=0$ with 
vector field $\alpha_0^*$, we have 
\begin{equation}
\frac{\partial\varepsilon}{\partial\alpha_0}+\frac1{\alpha_0} a_0\varepsilon
+A(x,\alpha_0)\varepsilon=0, \label{eq:3.1}
\end{equation} 
where $A(x,\alpha_0)=\sum_{H\in{\mathcal{A}}, H\not=H_0}
\frac{(\alpha_0^*,\alpha_H)}{\alpha_H}a_H$. 
Notice that 
\begin{align*}
s_0A(x\circ s_0^{-1},\alpha_0\circ s_0^{-1})s_0^{-1}
=&\sum_{H\not=H_0}\frac{(\alpha_0^*,\alpha_H)}{\alpha_H\circ s_0^{-1}}
s_0a_Hs_0^{-1}\\
=&\sum_{H\not=H_0}\frac{(s_0(\alpha_0^*),\alpha_{s_0(H)})}
{\alpha_{s_0(H)}}a_{s_0(H)}\\
=&\zeta_0 A(x,\alpha_0).
\end{align*}
Hence, 
\[B(x,\alpha_0):=\alpha_0 A(x,\alpha_0)\]
is $W_0$-invariant, i.e. 
$wB(x\circ w^{-1},\alpha_0\circ w^{-1})w^{-1}
=B(x,\alpha_0)$ for any $w\in W_0$, 
and satisfies $B(x,0)=0$. 
It follows that $B(x,\alpha_0)$ can be expressed as 
\begin{equation} 
B(x,\alpha_0)=\sum_{n=1}^\infty B_n(x)\alpha_0^n. \label{eq:3.2}
\end{equation} 
By definition, $s_0B_n(x)s_0^{-1}=\det^n(s_0)B_n(x)$. 

The equation \eqref{eq:3.1} is equivalent to 
\begin{equation} 
\alpha_0\frac{\partial\varepsilon}{\partial\alpha_0}
+a_0\varepsilon+B(x,\alpha_0)\varepsilon=0. \label{eq:3.3}
\end{equation} 

Let 
\begin{equation} 
\varepsilon(x,\alpha_0)=\alpha_0^c\sum_{n=0}^\infty\varepsilon_n(x)\alpha_0^n 
\label{eq:3.4}
\end{equation}
be a solution of \eqref{eq:3.3} satisfying $\varepsilon_0(x)=\varepsilon_0$. 
By \eqref{eq:3.2}, \eqref{eq:3.3} and \eqref{eq:3.4}, 
we have 
\begin{equation} 
\sum_{n=0}^\infty\left(\text{\hskip-2mm}\left(c+n+\sum_{l=0}^{e_0-1}k_{H_0,l}
\sum_{w\in W_0}{\det}^l(w)w\right)\varepsilon_n(x)
+\sum_{l=1}^nB_l(x)\varepsilon_{n-l}(x)\right)\alpha_0^{n+c}=0.
\label{eq:3.5}
\end{equation} 
By this equation, we have $c=-e_0k_{H_0,j}$ since 
$w\varepsilon_0=\det^{-j}(w)\varepsilon_0$ ($w\in W_0$). 

For a moment, assume that $e_0k_{H_0,l}\not\equiv e_0k_{H_0,l'}$ (mod $\Z$) 
for $l\not\equiv l'$ (mod $e_0\Z$). 
Then $n-e_0k_{H_0,j}+\sum_{l=0}^{e_0-1}k_{H_0,l}
\sum_{w\in W_0}\det^l(w)w$ is invertible 
since the eigenvalues of $w$ are $\det^{-l'}(w)$ 
($l'=0,\dots,e_0-1$). 
It follows that equation \eqref{eq:3.5} is uniquely solved and 
we have $w\varepsilon_n(x)=\det^{n-j}(w)\varepsilon_n(x)$ 
for $w\in W_0$ by induction. 
Hence 
\begin{equation}
\varepsilon(x,\alpha_0)=\alpha_0^{j-e_0k_{H_0,j}}\sum_{n=0}^\infty
\varepsilon_n(x)\alpha_0^{n-j}.
\label{eq:3.6} 
\end{equation} 
Since the series $\sum_{n=0}^\infty\varepsilon_n(x)\alpha_0^{n-j}$ is 
$W_0$-invariant, 
we have proved Theorem~\ref{thm:bmro} for generic $k$. 
By the continuity of $\varepsilon(x,\alpha_0)$ 
with respect to $k$, 
Theorem~\ref{thm:bmro} is proved. \qed
\end{pf}

\begin{defn}\label{defn:modi-tau-m}
{\rm Let \index{ei@$\varepsilon_i$} 
$\{\varepsilon_i; i=1,\dots,l=\deg\tau\}$ be a local basis of 
${\mathcal{E}}^{\nabla_\tau(k)}$ around $x_0\in X^\reg$ and 
let \index{sj@$\sigma_j$} 
$\{\sigma_j; j=1,\dots,l\}$ be a basis of 
$(\harm(\tau)\otimes E^*)^W$, as in the proof of Proposition
\ref{prop:nature of mintau}.  
For $k=(k_{C,j})_{C,j}$, we define a matrix \index{mk@$M(k)$} 
\[M(k)=(\sigma_j(\varepsilon_i))_{i,j=1,\dots,l}.\]}
\end{defn}

\begin{cor}\label{cor:tau-m} 
\begin{enumerate}
\item
The matrix \index{mzero@$M_0$} $M_0:=M(k=0)$ is a 
minimal $\tau$-matrix. 
\item
For all $k$, $M(k)$ is a nonsingular matrix. 
For any element $b$ of $B$, we define a matrix
\index{taukb@$\tau(k)(b)$} 
$\tau(k)(b)$ by 
\[\mu(b)M(k)=\tau(k)(b)M(k),\]
where \index{mub@$\mu(b)$} $\mu(b)M(k)$ is analytic continuation along
$b$. 
Then $\tau(k)$ is the monodromy representation of 
${\mathcal{E}}^{\nabla_\tau(k)}$ (cf. Theorem~\ref{thm:bmro}), written 
as matrix with respect to the 
local basis $\{\varepsilon_i; i=1,\dots,l=\deg\tau\}$. In 
particular, $\tau(k)(b)$ is an element of $GL(l,S)$ 
where $S$ denotes the ring of entire functions in $k$.
\item 
Let $\tau$ be irreducible.
The Euler vector field $\mathbb E$ acts on $M(k)$ in the following way:
\begin{equation}
{\mathbb E}M(k)=M(k)C_{\mathbb E}^M(k)
\label{eq:gegmin}
\end{equation}
with \index{cmek@$C_{\mathbb E}^M(k)$} 
\[C_{\mathbb E}^M(k)=C_{\mathbb E}^{M_0}-s(\tau,k)\id.\]
and where the scalar $s(\tau,k)$ is given by: 
\index{stauk@$s(\tau,k)$} 
\[s(\tau,k):=\frac1{\deg\tau}\sum_{C\in{\mathcal{C}}}
e_C|C|\sum_{j=0}^{e_C-1}k_{C,j}n_{C,j}^\tau.\]
\item
We have 
\[\det(M(k))=const.\prod_{C\in{\mathcal{C}}}\pi_C^{\epsilon_C(k)},\]
where \index{eck@$\epsilon_C(k)$} 
\[\epsilon_C(k)=\sum_{j=0}^{e_C-1}jn_{C,j}^\tau
-e_C\sum_{j=0}^{e_C-1}k_{C,j}n_{C,j}^\tau.\]
\item
If $k_{C,j}$ is an integer for every $C$ and $j$, then $\tau(k)(b)$ 
only depends on the image of $b$ in $W$ (cf. Theorem~\ref{thm:bmr}(iii)).  
This induces a map from $\Z^{\sum_{C\in{\mathcal{C}}}e_C}$ to 
the space $\funct(\rep(W),\rep(W))$ of functors from 
$\rep(W)$ to $\rep(W)$ by \index{gammak@$\gamma(k)$} 
\begin{align*}
\Z^{\sum_{C\in{\mathcal{C}}}e_C}\ni& k=(k_{C,j})\\
&\mapsto(\gamma(k): \tau\mapsto\tau(-k))
\in \funct(\rep(W),\rep(W)).
\end{align*}
Moreover, the matrix elements of $M(k)$ are rational functions on $V$ 
in this situation, with 
their poles contained in $\mathcal{A}$. 
\item
If $k$ is integral, 
then the local data (cf. \eqref{eq:rest}) $n_{C,j}^\tau$ and $n_{C,j}^{\tau(k)}$ 
coincide for any integral $k$. 
\end{enumerate}
\end{cor}
\begin{pf}
(i) This is a direct consequence of Definition \ref{defn:modi-tau-m}.
 
(ii) The regularity of $M(k)$ is a direct consequence of its definition. The 
remaining statement follows from Theorem \ref{thm:bmro}.

(iii) The Euler vector field $\mathbb E$ acts on ${\mathcal{E}}^{\nabla_\tau(k)}$. 
By definition, 
$\nabla_\tau(k)_{\mathbb E}={\mathbb E}+\langle{\mathbb E},\omega\rangle$ is $0$ on 
${\mathcal{E}}^{\nabla_\tau(k)}$. 
On the other hand, $\langle{\mathbb E},d\log\alpha_H\rangle=1$ 
for any $H\in{\mathcal{A}}$, 
and $\langle{\mathbb E},\omega\rangle=\sum_{H\in{\mathcal{A}}}a_H$ 
is an element of the 
center of $\C[W]$, hence it acts on $E$ by scalar multiplication since we 
assume $\tau$ to be irreducible in this item. 
This scalar is equal to $s(k,\tau)$,  
since 
\begin{align*}
\tr\sum_{H\in{\mathcal{A}}}a_H
=&\sum_{H\in{\mathcal{A}}}\sum_{j=0}^{e_C-1}k_{H,j}
\sum_{w\in W_H}{\det}^j(w)\tr_Ew\\
=&\sum_{H\in{\mathcal{A}}}\sum_{j=0}^{e_C-1}k_{H,j}
\sum_{w\in W_H}\sum_{j'=0}^{e_C-1}n_{C,j'}^\tau{\det}^{j-j'}(w)\\
=&\sum_{C\in{\mathcal{C}}}e_C|C|\sum_{j=0}^{e_C-1}k_{H,j}n_{C,j}^\tau.
\end{align*}
Therefore the Euler vector field $\mathbb E$ acts on ${\mathcal{E}}^{\nabla_\tau(k)}$ 
by the scalar 
$-s(\tau,k)$, 
and ${\mathbb E}M(k)=M(k)C_{\mathbb E}^M(k)$, as was claimed.

(iv) We may and will assume that $\tau$ is irreducible here. 
Since the series $\sum_{n=0}^\infty
\varepsilon_n(x)\alpha_{H_0}^{n-j}$ in \eqref{eq:3.6} is $W_{H_0}$-invariant, 
$\sigma(\varepsilon_n(x)\alpha_{H_0}^{n-j})$ is a $W_{H_0}$-invariant 
function for any $\sigma\in(\harm(\tau)\otimes E^*)^W$. 
It follows that this function is regular at $H_0$ because 
the pole of it is of order less than $e_{H_0}$. 

As in the proof of Proposition \ref{prop:nature of mintau} (iii), 
locally at $H\in C$, there exists a matrix $D\in GL(l,\C)$ 
such that 
\[DM(k)=M_H(k)M'(k),\]
where $M'(k)$ is a holomorphic $W_H$-invariant matrix and 
\index{mhk@$M_H(k)$} 
\begin{equation} 
M_H(k)=\diag(\alpha_H^{-k_{H,0}e_H}I_{n_{H,0}^\tau},
\alpha_H^{1-k_{H,1}e_H}I_{n_{H,1}^\tau},
\dots,\alpha_H^{e_H-1-k_{H,e_H-1}e_H}I_{n_{H,e_H-1}^\tau}).\label{eq:diag}
\end{equation}
It follows that 
\[\det M(k)=\prod_{C\in{\mathcal{C}}}\pi_C^{\epsilon_C(k)}\cdot r,\]
where, $r$ is a $W$-invariant, entire function. 
Using (iii) it follows that $\mathbb{E}r=0$ and (iv) is proved. 

(v) If $k_{C,j}$ is an integer for every $C$ and $j$, 
then $\mu(l_H)^{e_C}=\mu(l_H^{e_C})=1$ by Theorem~\ref{thm:bmro}. 
Hence $\tau(k)(b)$ only depends on the image of $b$ in $W$ by 
Theorem~\ref{thm:bmr}. This means that, with respect to the chosen 
basis of $E$, all the local sections of ${\mathcal{E}}^{\nabla_\tau(k)}$
extend to global holomorphic sections on $V^{\rm reg}$. 
Let us prove the rationality of $M(k)$ now. The connection  
$\nabla_\tau$ has simple poles at the reflection hyperplanes, and 
also at $\infty$. It follows that the pole orders at the reflection 
hyperplanes and at infinity are bounded by the eigenvalues of the 
residues of $\nabla_\tau$. Multiply a global flat section 
$\sigma$ with a suitable power of 
$\prod_{C\in{\mathcal{C}}}\pi_C$ so that the resulting product $\sigma^\prime$
is entire on $V$. Then the restriction of $\sigma^\prime$ to any complex line  
$L\subset V$ is a polynomial of degree $\leq N$ for some suitable $N\in \N$. 
Therefore, $\sigma^\prime$ is killed by all homogeneous constant coefficient 
differential operators of order $>N$, which implies that $\sigma^\prime$ is a 
polynomial. 
We conclude that $M(k)$ is a matrix of rational functions on $V$, 
whose poles are possibly at ${\mathcal{A}}$, and that $\tau(k)$ descends to 
a representation of $W$. 

We have 
\begin{align*}
\Z^{\sum_{C\in{\mathcal{C}}}e_C}\ni& k=(k_{C,j})\\
&\mapsto(\gamma(k): \tau\mapsto\tau(-k))
\in \funct(\rep(W),\rep(W)), 
\end{align*} 
since it is easy to verify the functoriality.

(vi) This follows from \eqref{eq:diag}. 
\qed
\end{pf}

\begin{rem} 
{\rm \begin{enumerate} 
\item
The functor $\gamma(k)$ does not respect $\otimes$, $\Hom$ 
and $*$ (contragredient). The reason for this is simply that 
the tensor product of two KZ connections is not equal to 
the KZ connection of the tensor product of the two representations 
of $W$ involved (and similarly for the other constructions from 
linear algebra that are mentioned).   
\item
If $W$ is a simply laced Coxeter group, 
the functor $\gamma(0,1)$ corresponds to the involution 
\index{i@$i$} \index{Lusztig's involution} 
$i: \hat{W}\rightarrow\hat{W}$ that was introduced by Lusztig, with the  
property that $i(\cdot)\otimes\det$ maps special representations to 
special representations. Similarly, for a Coxeter group that is not 
simply laced,  
the involution $i$ defined by Lusztig is equal to $\gamma((0,1),(0,1))$. 
 This involution $i$ plays a role in the study of 
cells and special representations of $W$.  
\end{enumerate}}
\end{rem}

\begin{conj}\label{ques}
For $k$ and $k'$ in $\Z^{\sum_{C\in{\mathcal{C}}}e_C}$, 
$\gamma(k)\gamma(k')$ equals 
$\gamma(k+k')$. 
\end{conj}

\begin{rem}\label{twijfels}
{\rm 
\begin{enumerate}
\item
Conjecture~\ref{ques} is clearly true if the following question, 
raised by Deligne-Mostow \cite{DM},  has an affirmative 
answer:
For a given irreducible representation $\tau$ of $W$, 
is it {\it uniquely} possible to deform $\tau$ to a 
family of representations $\tau_q$ of $B$ such that 
the eigenvalues of $\tau_q(l_H)$ satisfy the relation 
\[\prod_{j=0}^{e_C-1}(\tau_q(l_H)-\zeta_C^jq_{C,j})=0?\]
In turn, this question has an affirmative answer (locally in a 
neighbourhood of $q_{C,j}=1$ at least) if the 
so-called cyclotomic Hecke algebra of $W$, introduced by 
Brou\' e and Malle \cite{BM}, can be generated by 
$|W|$ elements over its coefficient ring. 
This was conjectured by Brou\'e, Malle and 
Rouquier \cite{BMR}, and checked in many cases. We will 
discuss these matters extensively in Section 6 and Section 7. 
\item
If Conjecture~\ref{ques} is true, then  
$\gamma(k)$ maps $\hat W$ to itself. This follows from that fact 
that $\gamma(k)$ respects direct sums and dimensions. If \ref{ques} 
is true, $\gamma(k)$ has inverse $\gamma(-k)$, 
but the   $\gamma(-k)$ pre-image of an 
irreducible representation must be irreducible by the above.   
In this situation we have an action $k\to\gamma(k)$ of the lattice 
$\Z^{\sum_{C\in{\mathcal{C}}}e_C}$ on $\hat W$. Let $I$ be the kernel 
of this action, and denote by \index{gw@$G_W$} 
$G_W$ the finite abelian group 
$\Z^{\sum_{C\in{\mathcal{C}}}e_C}/I$. We will identify this group 
with its image in ${\rm Per}(\hat W)$. 
In Section 6 we will show that this group is 
isomorphic with the geometric Galois group of the character field of the 
cyclotomic Hecke algebra, in the situation where the algebra is free 
over its coefficient ring. In particular we prove in this case that the character field 
is an {\it abelian} extension of the coefficient field of the cyclotomic 
Hecke algebra, assuming that the coefficient field contains $\C$.
\end{enumerate}}
\end{rem}

\section{Symmetries of fake degrees}

{\small In this section, we prove the ``fake degree symmetry'', 
which is the main theorem of this note.} 
\vskip3mm 

The key result is the following lemma, which shows how one 
may construct minimal matrices for $W$ using the KZ connection.
We use the notations introduced in Definition~\ref{defn:modi-tau-m}
and Corollary~\ref{cor:tau-m}.
\begin{lem}\label{lem}
Choose \index{b@$b=(b_C)_{C\in{\mathcal{C}}}$} 
$b=(b_C)_{C\in{\mathcal{C}}}$, where $b_C\in\{0,1,...,e_{C-1}\}$ for 
each $C\in {\mathcal{C}}$. 
Define \index{kb@$k_b$} 
$k_b\in \Z^{\sum_{C\in{\mathcal{C}}}e_C}$ by putting 
$(k_b)_{C,j}=1$ if $j\geq b_C$, and $(k_b)_{C,j}=0$ if $j< b_C$. 
Let \index{chib@$\chi_b$} 
$\chi_b$ be the one dimensional representation of $W$ 
associated with \index{pib@$\pi_b$} 
$\pi_b=\prod_{C\in{\mathcal{C}}}\pi_C^{e_C-b_C}$.
Then $\pi_b\cdot M(k_b)$ is a minimal matrix of type  
$\tau(k_b)\otimes\chi_b$ ($=\gamma(k_b)(\tau)\otimes\chi_b$).  
\end{lem}

\begin{pf}
By Corollary~\ref{cor:tau-m} (v), $M(k_b)$ is a 
rational matrix with poles in $\mathcal{A}$, 
and by \eqref{eq:diag} one sees that 
$\pi_b\cdot M(k_b)$ has polynomial entries. It is a 
nonsingular matrix of type 
$\tau(k_b)\otimes\chi_b$ by Corollary~\ref{cor:tau-m} (ii).  
Definition~\ref{defn:tau-m} (ii) is satisfied because of 
Corollary~\ref{cor:tau-m} (iii).  
It remains to prove the 
minimality. We will appeal to Remark~\ref{rem:mintau} (i) for the proof 
of minimality. By that remark it suffices to compute the degree of the 
determinant of $\pi_b\cdot M(k_b)$.
From Corollary~\ref{cor:tau-m} (iv) and (vi) we find that the degree equals 
$\sum_C|C|\epsilon^\prime_C(k_b)$ with 
\begin{align*}
\epsilon^\prime_C(k_b)
&=\epsilon_C(k_b)+(e_C-b_C){\rm dim}(\tau)\\
&=\sum_{j=0}^{e_C-1}jn_{C,j}^\tau-e_C\sum_{j=b_C}^{e_C-1}n_{C,j}^\tau+
(e_C-b_C){\rm dim}(\tau)\\
&=\sum_{j=0}^{b_C-1}(j+e_C-b_C)n_{C,j}^\tau+\sum_{j=b_C}^{e_C-1}(j-b_C)
n_{C,j}^\tau\\
&=\sum_{j=0}^{b_C-1}(j+e_C-b_C)n_{C,j}^{\tau(k_b)}+\sum_{j=b_C}^{e_C-1}(j-b_C)
n_{C,j}^{\tau(k_b)}\\
&=\sum_{j=0}^{e_C-1}jn_{C,j}^{\tau(k_b)\otimes\chi_b}
\end{align*}
This shows the minimality of $\pi_b\cdot M(k_b)$ by 
Remark~\ref{rem:mintau} (i).\qed
\end{pf}

\index{fake degree symmetry} 
\begin{thm}[Fake degree symmetry]\label{thm:main}
Let $\tau$ be irreducible.  
With notations as in Lemma~\ref{lem}, 
\[F_{\chi_b\otimes\tau(k_b)}(T)
=T^{N(\tau,b)}F_\tau(T).\]
where 
\[N(\tau,b)={\sum_{C\in{\mathcal{C}}}|C|\sum_{j=0}^
{b_C-1}(\frac{e_Cn_{C,j}^\tau}{\deg\tau}-1)}\]

\end{thm}
\begin{pf} 
This is a direct consequence of \eqref{eq:gegmin}, 
in view of Lemma~\ref{lem} and  Corollary~\ref{cor:2.6}.   \qed
\end{pf}

 

\section{Coxeter-like presentation of $W$} 

So far we used only one straightforward topological fact, the relation 
between the braid group $B$ and the pure braid group $P$, Theorem~\ref{thm:bmr}. 
We also avoided completely the use of the so called 
\index{cyclotomic Hecke algebra} cyclotomic Hecke algebra, although 
we already mentioned in Remark~\ref{twijfels} that this algebra plays an important role 
when one tries to say more about the meaning and properties of the functor $\gamma(k)$.
In the next sections we will discuss this in some detail. Unfortunately, one is 
forced to make the serious assumption that the cyclotomic Hecke algebra is a free 
algebra over its coefficient ring. This was conjectured by Brou\'e, Malle and 
Rouquier in \cite{BMR}, but at present no general proof is known. 
It has been proved for the infinite 
series by Ariki (see \cite{AK}) and by Brou\'e and Malle (see 
\cite{BM}) in about half of the exceptional cases. 

Let us now start the discussion of the 
cyclotomic Hecke algebra and its relation to the previous sections. The results 
and ideas in this and the next section are entirely due to Brou\'e, 
Malle and Rouquier. A good reference for this material is \cite{BMR}. 

To give the relation between the braid group and the cyclotomic Hecke algebra 
one first needs {\it Coxeter-like} presentations of $W$, that are well behaved 
with respect to the braid group $B$. 

\begin{defn} {\rm A \index{Coxeter-like presentation of $W$} 
Coxeter-like presentation of $W$ is a presentation of $W$  
given by a minimal set of generators \index{s@${\cal S}$} 
${\cal S}\subset W$ such that   
\begin{enumerate} 
\item
The set ${\cal S}$ consists of reflections of $W$, 
and all relations are generated by homogeneous 
\index{braid relation} {\it braid relations} 
and the \index{order relation} {\it order relations} 
(relations of the form $s^d=e$, for $s\in {\cal S}$). 
\item
There exists a choice of $l_H$ (for any $H$ with $s_H\in {\cal S}$) 
such that $\{l_H\}_{H, s_H\in {\cal S}}$ with just the homogeneous 
relations from (i) form a presentation of $B$. 
\end{enumerate}}
\end{defn} 

\begin{rem}{\rm  
\begin{enumerate} 
\item
It is well known that Coxeter groups have such presentations (see \cite{Br}, 
\cite{Del}). 
\item
Brou\'e-Malle-Rouquier have shown that there exists such a presentation 
of $W$ for all but $G_{24}, G_{27}, G_{29}, G_{31}, G_{33}$ 
and $G_{34}$. For these groups, the existence of such a presentation 
is conjectural. 
\end{enumerate}}
\end{rem}



\section{Cyclotomic Hecke algebras and monodromy} 

{\small In this section, we define the Hecke algebra $H_u(W)$ 
of finite complex reflection groups which have a Coxeter-like presentation.} 
\vskip3mm

Let \index{s@$S$} $S$ be the ring of entire functions in 
$(k_{C,j})_{C\in{\mathcal{C}}, j=0,\dots,e_C-1}$ 
and let \index{k@$K$} $K$ be its field of fractions. 
We also define \index{r@$R$} \index{ucj@$u_{C,j}^\pm$} 
$R=\Z[u_{C,j}^\pm; C\in{\mathcal{C}}, j=0,\dots,e_C-1]$ 
($u_{C,j}$ are indeterminant) and let \index{q@$Q$} 
$Q$ be its field of fractions. 
The following definition is from \cite{BMR}.
 
\begin{defn} 
{\rm The \index{topological Hecke algebra} 
(topological) Hecke algebra \index{huw@$H_u(W)$} $H_u(W)$ 
is the $R$-algebra $R[B]/J$, 
where $J$ is the ideal generated by the elements 
$\prod_{j=0}^{e_H-1}(l_H-u_{H,j})$.} 
\end{defn}

\begin{cor}\label{pres} 
If $W$ has Coxeter-like presentation, then the image 
\index{ts@$T_s$} 
$\{T_s\mid s\in {\cal S}\}$ of $\{l_H\mid  H\in{\mathcal{A}} 
{\text{ such that}} \enskip s_H=s\in {\cal S}\}$ 
generates $H_u(W)$ and is subject to the 
relations $\prod_{j=0}^{e_C-1}(T_s-u_{s,j})=0$. Together with 
the braid relations (which already hold in $B$) this   
is a presentation of $H_u(W)$. 
(Here $u_{s,j}=u_{H,j}$ if $s=s_H$.) 
\end{cor} 

With this definition, 
the next theorem follows directly from Theorem~\ref{thm:bmro}, except 
for the last assertion (iii).

\begin{thm}\label{thm:6.3}
We take Knizhnik-Zamolodchikov connection with values in  
$E=\C[W]$, and let \index{muk@$\mu(k)$} $\mu(k)$ be 
its monodromy representation. 
\begin{enumerate} 
\item 
For $b\in B$, $\mu(k)(b)$ coincides with the left multiplication 
by an element $\mu(k)_b\in(\C[W])^\times$. 
\item
Let us embed $R$ in $S$ by means of 
the substitution 
\[u_{H,j}=\det(s_H)^j\exp(-2\pi\sqrt{-1}k_{C,j}).\]
Then the representation $\mu(k)$ 
factorizes as follows: 
\[\matrix
S[B]&&\overset{\mu(k)}
\longrightarrow&&S[W]\\
&\searrow&&\nearrow&\\
&&S\otimes_R H_u(W)&&
\endmatrix.\]
A similar factorization holds for the matrix representations $\tau(k)$ defined 
in Corollary~\ref{cor:tau-m}(ii). These representations can be obtained by 
restriction of $\mu(k)$ to a minimal left ideal of $\C [W]$ of type $\tau$.
\item 
(Brou\'e-Malle-Rouquier \cite{BMR}) 
If $W$ has a Coxeter-like presentation and $H_u(W)$ can be generated by $|W|$ elements 
as $R$-module, then $\mu(k)$ induces an $K$-algebra isomorphism 
\begin{equation}\label{ss}
K\otimes_RH_u(W)\overset{\sim}\rightarrow K[W].
\end{equation}
Moreover, in this situation $H_u(W)$ is free over $R$ of rank $|W|$.
\end{enumerate} 
\end{thm} 

\begin{rem}\label{rem:core}{\rm  
The assumption of Theorem~\ref{thm:6.3} (iii) is conjectured to be true in all cases and 
is known for the basic series by the work of Ariki and Koike \cite{AKo} 
and Ariki \cite{AK}, and the work of 
Brou\'e and Malle \cite{BM} in many of the exceptional cases.  
}
\end{rem} 

\begin{hyp}\label{hyp} 
{\rm  From now on we will assume that $W$ has a Coxeter-like
presentation and that $H_u(W)$ can be generated by $|W|$ elements, as required 
in Theorem~\ref{thm:6.3} (iii). 
}
\end{hyp}

\begin{cor}\label{cor:split} The algebra 
$Q\otimes H_u(W)$ is a semisimple algebra of rank $|W|$, which splits over 
$K\supset Q$.
\end{cor}

Let \index{f@$F$} $F$ denote the character field of $Q\otimes H_u(W)$, 
which we define as the 
subfield of $K$ generated by the values of the irreducible characters of 
$K\otimes H_u(W)$ on $Q\otimes H_u(W)$. It is a subfield of $K$ containing 
$Q$. Note that the characters take values in $S$ on an $R$-basis of 
$H_u(W)$, by Theorem~\ref{thm:bmro} and Theorem~\ref{thm:6.3}(ii), (iii). 
In fact, the values of these characters on $H_u(W)$ are even in the integral closure of  
$R$ in $S$, by a well known argument due to Steinberg 
(for example, see \cite{C}, Proposition 10.11.4).
In particular, $F$ is a Galois extension of $Q$.
 Let \index{gf@$G_F$} $G_F$ denote its Galois group.
This Galois group acts naturally on the irreducible characters of $K\otimes H_u(W)$. 
The evaluation homomorphism $f\to f(0)$ of $S$,  
applied to the characters, defines a bijection between 
the irreducible characters of $K\otimes H_u(W)$ and of $W$ (by Tits' specialization 
theorem). Hence we arrive at a certain action of $G_F$ on $\hat W$. The Galois 
group of $\C\cdot F$ over $\C\cdot Q$ is called 
the \index{geometric Galois group} geometric Galois group, 
which we will denote by \index{gfgeom@$G_F^{geom}$} $G_F^{geom}$. 
By restriction to $F\subset \C\cdot F$ we 
have a canonical homomorphism $G_F^{geom}\to G_F$, 
which is obviously injective. 
We identify $G_F^{geom}$ with its image in $G_F$. 
Via this embedding of $G_F^{geom}$ in $G_F$ we 
have now defined an action of $G_F^{geom}$ on $\hat W$, 
which we will denote by \index{al@$\alpha$} $\alpha$. 
Note that $\alpha$ is an injective homomorphism $G_F^{geom}\to{\rm Per}(\hat W)$.

\begin{thm}\label{thm:pol} 
If Hypothesis~\ref{hyp} holds, then the field $\C\cdot F$ is contained in a field that 
is obtained from $\C\cdot Q$ by adjoining  radicals of the form $(u_{C,j})^{1/N}$,
for a suitable $N$. In particular,
$\C\cdot F$ is an abelian extension of $\C\cdot Q$.
\end{thm} 
\begin{pf} 
As was explained above, the characters of $K\otimes H_u(W)$ take values in the 
integral closure of $R$ in $S$ on $H_u(W)$.     
The lattice \index{l@$L$} 
$L=\Z^{\sum_{C\in{\mathcal{C}}}e_C}$ (as in Corollary~\ref{cor:tau-m} (v)) 
acts on $S$ by means of the action \index{bl@$\beta(l)$} 
$(\beta(l)f)(k)=f(k-l)$. Notice that 
$Q$ is fixed for this action!
By Theorem~\ref{thm:6.3}(ii),  for every $\tau\in \hat W$ we have 
a representation $\tau(k): S\otimes H_u(W)\to {\rm Mat}(d_\tau\times d_\tau, S)$ 
(where \index{dt@$d_\tau$} $d_\tau$ is the degree of $\tau$)  
and these representations exhaust the equivalence classes of irreducible 
representations of $K\otimes H_u(W)$ (by Theorem~\ref{thm:6.3}(iii)). 
We define  $\beta(l)$ $(l\in L)$ 
in the obvious way on $S\otimes H_u(W)$ and on ${\rm Mat}(d_\tau\times d_\tau, S)$ 
(i.e. coefficient-wise) and then define $\beta(l)(\tau)=\beta(l)\circ\tau\circ\beta(-l)$. 
Obviously, $\beta(l)(\tau)$ is $S$-linear, and it is also a representation 
of $S\otimes H_u(W)$ because $\beta(l)$ is a ring homomorphism of $S\otimes H_u(W)$. 
Using Tits' specialization theorem, applied 
to specializations of $S$ at lattice points, we find that the new representation so 
obtained by translation, is irreducible. In this way, we have an action $\beta$ of the 
lattice $L$ on the finite set ${\rm Irr}(K\otimes H_u(W))$. Hence there exists 
a $N\in\N$ such that $N\cdot L$ fixes every representation. This means that a character value 
$\chi_\tau(b)$ is periodic with respect to the lattice $N\cdot L$. Therefore it can be 
viewed as a univalued function of the coordinates
\index{zcj@$z_{C,j}$} $z_{C,j}=(u_{C,j})^{1/N}$, which is 
as such still a holomorphic function outside the coordinate hyperplanes 
$z_{C,j}=0$, 
and integral over $Q$, hence certainly integral 
over the larger ring of polynomials in $z_{C,j}$. But  
this polynomial ring is integrally closed in the field of meromorphic functions
in the $z_{C,j}$, 
because a solution of a monic equation is bounded in norm by the 
sum of the norms of the  coefficients of the equation (including the 
top coefficient 1).  
Our conclusion is that  $\chi_\tau(b)$ is a polynomial in the 
coordinates $z_{C,j}=(u_{C,j})^{1/N}$. \qed
\end{pf}

\begin{cor}
\begin{enumerate}
\item Hypothesis~\ref{hyp} implies that Conjecture~\ref{ques} holds. In particular, 
as was mentioned in Remark~\ref{twijfels} (ii), the functors $\gamma(l)$ of 
Corollary~\ref{cor:tau-m}(v) define an action of $L$ on $\hat W$.
\item The action $\alpha$ of $G_F^{geom}$ on $\hat W$ defines an isomorphism of 
$G_F^{geom}$ with the group $G_W$ described in Remark~\ref{twijfels} (ii).  
\end{enumerate} 
\end{cor}
\begin{pf}
In fact, everything is clear by the remark that the action $\beta$ of $L$ defined in 
the proof of Theorem~\ref{thm:pol} is equal to the action $\gamma$ when we 
identify ${\rm Irr}(K\otimes H_u(W))$ and $\hat W$ via the specialization $f\to f(0)$ 
of $S$, as always. From Theorem~\ref{thm:pol} is is clear that $G_F^{geom}$ is 
the quotient of $L/(N\cdot L)$ by the subgroup that fixes all the 
irreducible characters via the action $\beta$, and this is isomorphic to $G_W$ 
by the definition of $G_W$.\qed
\end{pf}



\section{Applications}

The phenomenology in the field of \index{cyclotomic Hecke algebra} 
cyclotomic Hecke algebras is much 
further developed than the theory, and the present paper does not 
change this situation very much! Recently, Gunter Malle made a thorough 
study of the character and splitting fields of representations of cyclotomic 
Hecke algebras in \cite{M}. His study is based on the theory developed in 
\cite{AK} for the 
case of the cyclotomic Hecke algebra of the groups $G(m,p,n)$, and on a  
case by case analysis of the primitive groups (always assuming 
Hypothesis~\ref{hyp} of course). There are many intriguing observations 
in Malle's paper \cite{M}, and some of these are closely related to the 
questions we have discussed here.  {\it One remark has to be made beforehand. 
In the comparison of results of \cite{M} and the results of the present paper 
one must realize that we are discussing the topological Hecke 
algebra here, with its presentation as described in Corollary~\ref{pres}.  
However, Malle uses the \index{abstract Hecke algebra} 
abstract Hecke algebra, defined with a similar presentation 
but where the topological braid group $B$ is replaced by the abstract
braid group 
associated with a certain choice of a presentation of $W$ (by removing the 
order relations from the presentation of $W$). In other words, only when the 
presentation of $W$ that Malle uses is a Coxeter-like presentation  it 
is clear that we are discussing the same algebra. Recall that there are still 
some groups for which the existence of a Coxeter-like presentation is
not known. See section 5.} With this understood, let me list the main 
facts revealed in Malle's paper, and comment on these from the point 
of view of the theory in this paper.

\begin{enumerate}
\item The character field $F$ is a regular, abelian extension of 
$k(u_{C,j})$, 
where \index{k@$k$} $k$ denotes the splitting field of the group
algebra $\Q[W]$. 
(It is known 
that $k$ is the character field of the reflection representation of $W$, 
by a result of \cite{Bs}). 
This means that $G_F^{geom}={\rm Gal}(F/k(u_{C,j}))$, 
and that this group is abelian.  The equality of these 
two Galois groups is based on the (empirical!) fact that 
$k$ contains all the roots of unity of order $d$,  
whenever $d$ is the order of an element of $G_F^{geom}$ ({\it loc. cit.} Cor. 4.8).
\item The order of $G_F^{geom}$ can be arbitrarily large (but of course subject 
to the condition imposed by (i)), even if all the reflections 
in $W$ have order 2 (cf {\it loc. cit.} Ex. 4.5). The orders of elements of 
$G_F^{geom}$ do not  
necessarily divide the order of the center $Z(W)$ (but in the ``well-generated''
case they do, see {\it loc. cit.} Prop. 7.2). 
\item A Beynon-Lusztig type of ``semi-palindromicity'' for fake
  degrees of rational representations of the 1-parameter cyclotomic
  Hecke algebra ({\it loc. cit.} Thm. 6.5). 
\item The group $W$ can be generated by $n$ complex reflections if and
  only if the fake degree of the reflection representation is
  semi-palindromic ({\it loc. cit.} Prop. 6.12). 
\end{enumerate}

We have not much to say about (i), since the methods used here do not seem fit 
for the study of the full Galois group of $F$. The only 
thing we have proved is the fact that $G_F^{geom}$ is 
abelian (Theorem~\ref{thm:pol}). 

It is implied by (i) that the order of 
every element of $G_F^{geom}$ is a divisor of the order of the group of roots 
of unity in $k$. Looking at the explicit results of Malle and of 
Brou\' e and Malle \cite{BM} in the primitive cases, and at (ii) mentioned 
above, this seems to be the only general rule. 
From the point of view of the present paper this fact is rather 
mysterious, except in some cases where the order of all nontrivial 
elements in 
$G_F^{geom}$ is two. This includes all the Coxeter cases. 
The case of Coxeter groups was dealt with in 
the paper \cite{O}, but since the argument given there is incomplete, 
we will fill in the details here. 

We start with some elementary relations:

\begin{prop}\label{order} 
Let \index{gcj@$g_{C,j}$} 
$g_{C,j}\in G_W=G^{geom}_F$ be the image under $\gamma$ of 
the basis element \index{bcj@$b_{C,j}$} $b_{C,j}$ of the lattice $L$, 
defined by 
$(b_{C,j})_{C',j'}=-\delta_{C,C'}\delta_{j,j'}$. 
Let \index{mc@$m_C$} $m_C$ denote the permutation 
of irreducible representations of $W$ given by:
$m_C(\tau)=\tau\otimes\chi_C$, where 
\index{chic@$\chi_C$} $\chi_C$ is the linear character associated to
the pseudo invariant $\pi_C$. 
We have the following relations:
\begin{enumerate}
\item The $g_{C,j}$ mutually commute.
\item $g_{C,0}\dots g_{C,e_C-1}={\rm id}$.
\item $m_C^{e_C}={\rm id}$.
\item $g_{C,j}m_{C'}=m_{C'}g_{C,j}$ ($C\not=C'$).
\item $g_{C,j}m_C=m_Cg_{C,j+1}$ (cyclic).
\end{enumerate}
In particular, we have ($1\leq j\leq e_C$): $(m_Cg_{C,1}\cdots g_{C,e_C-1})^j=
m_C^jg_{C,j}\cdots g_{C,e_C-1}$. Thus the order of $m_Cg_{C,1}\cdots g_{C,e_C-1}$
divides $e_C$. (Note that the powers of  $m_Cg_{C,1}\cdots g_{C,e_C-1}$ are the 
operations on the representation $\tau$ that occur in the fake degree symmetry 
formula Theorem~\ref{thm:main}.)
\end{prop}
\begin{pf} In the notation of Definition~\ref{kztau}, 
assertion (ii) follows from the remark (verified by 
means of direct computation) that the flat sections of 
$\nabla_\tau(b_{C,0}+\dots +b_{C,e_C-1})$ are of the form 
$\pi_C^{e_C}\cdot e$ with $e\in E$. This has type $\tau$ since 
$\pi_C^{e_C}$ is $W$-invariant. Likewise, (iv) and (v) follow 
from a direct computation of the tensor product of the 
KZ-connection $\nabla_\tau(k)$ and the (trivial) KZ-connection 
$\nabla_{\chi_C}(0)$. It is left to the reader. The remaining 
statements are trivial.\qed
\end{pf}

\begin{rem} 
{\rm 
It is informative to have a look at Malle's table 7.1, 
displaying all the 
irrationalities for the primitive groups, at this point. For example, one 
can check Proposition~\ref{order}(ii) directly from this table: when we 
specialize $u_{C,j}\to u_C$ everything becomes rational in the $u_C$. 
}
\end{rem}

Now we may reprove the following well known results:

\begin{prop}\label{ouwekoek} Let $W$ be a Coxeter group, 
or of type $G(4,2,n)$ ($n>2$), or $G_{31}$, in which case we need 
hypothesis~\ref{hyp}, as always. 
Then the elements $g_{C,j}$ ($j=0,1$) have order 
2, and $g_{C,0}=g_{C,1}$. 
(In other words, $\C \cdot F$ is a subfield of the field obtained from 
$\C\cdot Q$ by adjunction of the square roots $\sqrt{q_{C,0}q_{C,1}}$.)    
\end{prop} 
\begin{pf} The point is that there exists a 
\index{Kazhdan-Lusztig involution} Kazhdan-Lusztig involution 
\index{j@$j$} $j$ of 
the (cyclotomic) Hecke algebra, defined by $j(T_i)=-q_{C,0}q_{C,1}T_{i'}^{-1}$ 
where $s_i$ and $s_{i'}$ belong to hyperplanes in $C$. In fact, if $W$ is 
Coxeter we may take $i=i'$, but for the complex cases 
we have to flip elements 
according to the symmetry axes in the diagram, in order to 
preserve the circular relations. This induces an involution on the 
irreducible representations, also denoted by $j$, defined by 
$j(\tau)=\tau\circ j$.
It is simple to see that $j(\tau)=\tau\otimes {\rm det}=\prod_Cm_C(\tau)$ 
in the Coxeter case. In the complex cases we mentioned this is also 
true because 
the aforementioned flip of generators is an {\it inner} automorphism.  
(For $G_{31}$ this follows from a calculation in \cite{BM}, Bemerkung 6.5, 
and for the infinite series $G(4,2,n)$ it is similar.) 

On the other hand, $j$ obviously commutes with the elements $g_{C,i}$ 
since it is rational in the $q_{C,i}$. Combining this with 
Proposition~\ref{order} we obtain the result.\qed
\end{pf}

There is a bewildering number of groups acting on $\hat W$: $Gal(F/Q)$, the 
group of linear characters (via tensoring), complex conjugation, 
diagram automorphisms, and the symmetric groups \index{sc@$\S_C$} $\S_C$ of 
permutations of $\{u_{C,j}\mid 0\leq j\leq e_C-1\}$. There are 
some easy observations about the groups that some of these actions 
generate.  
For example, the actions of symmetric group $\S_C$ and the lattice 
\index{lc@$L_C$} $L_C$ 
(i.e. all the parameters are 0 for hyperplanes not in $C$) combine to give 
an action of the affine Weyl group \index{scaff@$\S_C^{\rm aff}$} 
$\S_C^{\rm aff}$. Or, as remarked in 
\cite{M}, complex conjugation and $g_{C,0}$ generate an action 
of the infinite dihedral group. Another example is the cyclic group of order 
$e_C$ described in Proposition~\ref{order}. Maybe it is important 
to investigate this systematically, but I did not see anything 
useful other than Proposition~\ref{ouwekoek}

The next issue (iii) from Malle's paper, the Beynon-Lusztig type of 
``semi-palindromicity'' of the fake degrees of rational representations, 
can be fully understood in terms of our results. The Galois operation 
\index{del@$\delta$} $\delta$ he introduces is easily identified in
our notations as 
$\delta=(\prod_Cg_{C,0})^{-1}$ (recall how $R$ is 
embedded in $S$, see Theorem~\ref{thm:6.3}). 
With this notation, Malle observes the 
following:
\begin{prop}(\cite{M}, Theorem 6.5)\label{pali}
Let \index{r@$\mathcal R$} $\mathcal R$ denote the set of reflections
in $W$. 
In the notations of Definition~\ref{fake}, we have
\begin{equation}
R_\tau(T)=T^cR_{\delta(\overline{\tau})}(T^{-1})
\label{pal}
\end{equation}
where $c=\#\mathcal R-\sum_{r\in\mathcal R}\chi_\tau(r)/\chi_\tau(1)$.
\end{prop}
\begin{pf} 
Replace $T$ by $T^{-1}$ in \eqref{pal}, multiply 
both sides by $T^{\#\mathcal R}$ and rewrite the left 
hand side using the well known and elementary formula:
\[
T^{\#\mathcal R}R_{\tau}(T^{-1})=F_{\tau\otimes{\rm det}}(T).
\]
On the right hand side we use Definition~\ref{fake} and 
the observation $\overline{\delta(\overline{\tau})}=
\delta^{-1}(\tau)$, to obtain
\[
F_{\tau\otimes{\rm det}}(T)=
T^{N-c}F_{\delta^{-1}(\tau)}(T)
\]
Finally replace $\tau$ by $\delta(\tau)$ and use 
Proposition~\ref{order}(ii) to see that (in the notation 
of Theorem~\ref{thm:main})
$\delta(\tau)\otimes{\rm det}=\tau(k_b)\otimes \chi_b$
where the vector $b$ is defined by 
$b_C=1\forall C$. The result now follows from Theorem~\ref{thm:main}.\qed
\end{pf}  
 
Finally, we have quoted assertion (iv) of Malle's paper because it is such  
a nice result related to fake degees. 
It goes back to Orlik and Solomon, who proved it by 
inspection. Malle suggests an alternative 
approach to this result, based on two empirical facts.  
First, the reflection character of the ``1-parameter 
specialisation'' of the Hecke algebra 
is rational if and only if $W$ is generated by $n$ reflections. 
Next, the reflection character of the ``1-parameter 
specialisation'' of the  Hecke algebra is rational if and 
only if its fake degree is semi-palindromic 
(see \cite{M}, Corollary 4.9 and Proposition 6.12). 
Hence \ref{pali} gives an a priori proof of one out of the 
four implications in these two statements.

\newpage
\printindex

\end{document}